\documentclass[fleqn,10pt,twocolumn]{SICE18}  

\newcommand{\Sym}{\operatorname{Sym}}
\newcommand{\Skew}{\operatorname{Skew}}

\newcommand{\SO}{\operatorname{SO}(3)}
\newcommand{\MI}{{\mathbb I}}
\newcommand{\dimu}{{k}}
\newcommand{\coV}{{k_{\rm e}}}

\title{
Observer-Based Controller Design for Systems on Manifolds in  Euclidean Space
}

\author{Dong Eui Chang
}

\speaker{Dong Eui Chang}

\affils{School of Electrical Engineering, Korea Advanced Institute of Science and Technology, Daejeon, Korea\\
(E-mail: dechang@kaist.ac.kr)
}

\abstract{
A method of designing observers and observer-based tracking controllers is proposed for  nonlinear systems on manifolds  via  embedding into Euclidean space and transversal stabilization. Given a system on a manifold, we first embed the manifold and the system into Euclidean space and extend the system dynamics to the ambient Euclidean space in such a way that the manifold becomes an invariant attractor of the extended system, thus securing the transversal stability of the manifold in the extended dynamics. After the embedding, we design state observers and observer-based controllers for the extended system in one single global  coordinate system in the ambient Euclidean space, and then restrict them to the original state-space manifold to produce observers and observer-based controllers for the original system on the manifold. This procedure has the merit that  any existing control method that has been developed in Euclidean space can be applied globally to systems defined on nonlinear manifolds, thus making nonlinear controller design on manifolds easier. The detail of the method is demonstrated on the fully actuated rigid body system.}


\begin{document}

\maketitle

\section{Introduction}

We have studied in \cite{Ch17} about tracking controller synthesis for systems defined on manifolds via  embedding into Euclidean space, transversal stabilization, and linearization.  The main idea of this method is as follows. Given a control system on a manifold $M$, first embed $M$ into  Euclidean space ${\mathbb R}^n$ and then extend the control system dynamics from $M$ to ${\mathbb R}^n$ in such a way that $M$ becomes an invariant attractor of the extended or ambient system in ${\mathbb R}^n$, thus securing the transversal stability of $M$ in the extended dynamics. As a result, we can conveniently design controllers using one single global Cartesian coordinate system in ${\mathbb R}^n$ for the original system on $M$. In particular, the linearization of the system dynamics along a reference trajectory on $M$ can be carried out {\it globally}, without relying on multiple local charts,  in a Cartesian coordinate system in ${\mathbb R}^n$ to easily design tracking controllers for the system on $M$. This methodology for controller synthesis is well illustrated with the fully actuated rigid body system in \cite{Ch17}. We refer the reader to \cite{Ch11,ChJiPe16} for other applications of the embedding technique in the context of optimal control and geometric numerical integrations.

In this paper, we continue the program of embedding into Euclidean space, transversal stabilization, and linearization in order to  construct state observers and observer-based tracking controllers  for systems defined on manifolds. Given a system on a manifold $M$, we embed it into  Euclidean space ${\mathbb R}^n$ in such a way that $M$ becomes an invariant attractor of the embedded system in ${\mathbb R}^n$, and then design both controllers and observers for the system through linearization along a given reference trajectory on $M$ in one single global Cartesian coordinate system in ${\mathbb R}^n$ after which these observers and controllers are combined to produce observer-based controllers for the original system on $M$; any other observer design technique in ${\mathbb R}^n$ can be utilized  although linearization is employed in this paper for its simplicity but high utility.  Thanks to the use of one single global Cartesian coordinate system, there is no need of change of coordinates along the reference trajectory. Also, all the conditions for (local) exponential stability of the linear time varying (LTV) observer laid out in \cite{Ka60,KaBu61} can be verified  {\it a priori} globally along the reference trajectory on $M$.

This paper is organized as follows. The program of  embedding into Euclidean space, transversal stabilization, and linearization is first reviewed and then the method for designing observers and observer-based tracking controllers is developed in the context of this program. Two kinds of Kalman-type observers are constructed:  one is LTV observers that estimate the tracking error and the other is high-gain observers that estimate the state of the system, both of which utilize Kalman's theory of linear systems \cite{Ka60,KaBu61}. The observer-based tracking controller design procedure is illustrated with the fully actuated rigid body system and a good tracking and state observation performance of the controller for  large initial tracking/observation errors is demonstrated with a simulation.  In addition to the two Kalman-type observers, a non-Kalman-type  observer is  constructed for the rigid body system.  For the sake of completeness of the paper and convenience to the reader, the theory of LTV systems is briefly reviewed in a modified but improved form in the Appendix. 

\section{Main Results}

\subsection{Review}
\subsubsection{Embedding and Transversal Stabilization}
We here review  the thechnique of embedding and transversal stabilization from \cite{Ch17}. Let $M$ be a regular manifold in some ${\mathbb R}^n$. Consider a control system $\Sigma_M$ on $M$ given by
\begin{equation}\label{Sigma:M}
\Sigma_M: \quad \left \{ 
\begin{array}{l}
 \dot x = X(x,u), \quad x \in M, u \in {\mathbb R}^\dimu,\\
 y= h(x)  \in {\mathbb R}^p.
\end{array}
\right .
\end{equation}
where $x$ is the state, $u$ the control,  and $y$ the output of the system.  Here we assume that the function $h$ is defined on ${\mathbb R}^n$. It is understood that $X(x,u) \in T_xM$ for all $(x,u) \in M\times {\mathbb R}^\dimu$.
Suppose that there is a control system  $\Sigma_{{\mathbb R}^n}$ in ${\mathbb R}^n$  given by
\begin{equation}
\Sigma_{{\mathbb R}^n}: \quad \left \{ 
\begin{array}{l}
 \dot x = X_e(x,u), \quad x \in {\mathbb R}^n, u \in {\mathbb R}^\dimu,\\
 y= h(x)  \in {\mathbb R}^p,
\end{array}
\right .
\end{equation}
where it is assumed that 
\[
X_e(x,u) = X(x,u) \quad \forall x\in M, u \in{\mathbb R}^\dimu,
\]
which means that $\Sigma_{{\mathbb R}^n}$ is an extension of $\Sigma_M$ to ${\mathbb R}^n$, and $\Sigma_M$ is a restriction of $\Sigma_{{\mathbb R}^n}$ to $M$.

Suppose that there is a non-negative function $\tilde V$ on ${\mathbb R}^n$ such that $M = \tilde V^{-1}(0)$ and
\begin{equation}\label{tilde:V:X}
\nabla \tilde V (x) \cdot X_e (x,u) =0 
\end{equation}
for all $(x,u) \in {\mathbb R}^n \times {\mathbb R}^\dimu$.
Consider the following system 
\begin{align}\label{Sigma:tilde}
\tilde \Sigma_{{\mathbb R}^n}: \quad \left \{ 
\begin{array}{l}
 \dot x = \tilde X_e(x,u), \quad x \in {\mathbb R}^{n}, u \in {\mathbb R}^\dimu,\\
 y= h(x)  \in {\mathbb R}^p,
\end{array}
\right .
\end{align}
where
\[
\tilde X_e(x,u) :=X_e(x,u) -\nabla \tilde V(x) \quad \forall x\in {\mathbb R}^n, u \in{\mathbb R}^\dimu.
\]
Since the minimum value $0$ of $\tilde V$ is attained on $M$, the gradient $\nabla \tilde V$ identically vanishes on $M$, which implies that $\tilde \Sigma_{{\mathbb R}^n}$ coincides with $\Sigma_M$ on $M$,  and that $M$ is an invariant manifold of the new system $\tilde \Sigma_{{\mathbb R}^n}$. A sharp distinction between the two systems $\tilde \Sigma_{{\mathbb R}^n}$ and $\Sigma_{{\mathbb R}^n}$ is that $M$ is a stable  invariant manifold of $\tilde \Sigma_{{\mathbb R}^n}$ due to the added term $-\nabla \tilde V$. The following new theorem provides a sufficient condition for asymptotic stability of $M$ in the transversal direction for the system $\tilde \Sigma_{{\mathbb R}^n}$. 
\begin{theorem}\label{theorem:attractive:M}
If there are positive numbers $b$  and $r$ such that $
b  \tilde V (x) \leq  \| \nabla \tilde V(x) \|^2 $
for all $x \in \tilde V^{-1}([0,r))\subset {\mathbb R}^{n}$, then $\tilde V^{-1}([0,r))$ is positively invariant for $\tilde \Sigma_{{\mathbb R}^n}$ and every flow of $\tilde \Sigma_{{\mathbb R}^n}$ starting  in $\tilde V^{-1}([0,r))$ converges to $M$ as $t\rightarrow \infty$. In particular, $\tilde V(x(t))  \leq \tilde V(x(0)) e^{-b t} $ for all $t\geq 0$ and $x(0) \in \tilde V^{-1}([0,r))$.
\end{theorem}

\subsubsection{Tracking Controller Design via Linearization in Ambient Euclidean Space}
We review from \cite{Ch17} the technique of tracking controller design via linearization in ambient Euclidean space after embedding. 
Consider a reference trajectory $x_0: [0,\infty) \rightarrow M$ for the system $ \Sigma_M$ on $M$ driven by a control signal $u_0: [0,\infty) \rightarrow {\mathbb R}^\dimu$, so that 
\[
\dot x_0(t) =  X(x_0(t), u_0(t)) \quad \forall t\geq 0.
\]
Our goal is to construct tracking controllers for this trajectory for the system $\Sigma_M$. 
By the construction of $\tilde \Sigma_{{\mathbb R}^n}$, this reference trajectory also satisfies the dynamics of $\tilde \Sigma_{{\mathbb R}^n}$, i.e.
\[
\dot x_0(t) =  \tilde X_e(x_0(t), u_0(t)) \quad \forall t\geq 0.
\]
Hence, we take the  strategy to design  tracking controllers for the ambient system $\tilde \Sigma_{{\mathbb R}^n}$ in ${\mathbb R}^n$ and then restrict them to $M$ to come up with  tracking controllers for  $\Sigma_{M}$.
For convenience, assume that both $x_0(t)$ and $u_0(t)$ are bounded signals.
Let 
\begin{align}
\Delta x &= x - x_0(t) \in {\mathbb R}^n,\label{def:Delta:x}\\
 \Delta u &= u - u_0(t) \in {\mathbb R}^\dimu, \label{def:Delta:u}\\
\Delta y &= y - h(x_0(t)) \in {\mathbb R}^p. \label{def:Delta:y}
\end{align}
Then, they satisfy
\begin{subequations}\label{tracking:dyn}
\begin{align}
\Delta \dot x &= \tilde X(x_0(t) + \Delta x,u_0(t) + \Delta u ) -\tilde X(x_0(t),u_0(t)),\label{tracking:dyn:a}\\
\Delta y &=  h( x_0(t) + \Delta x) - h(x_0(t)) \label{tracking:dyn:b}
\end{align}
\end{subequations}
or
\begin{subequations}\label{tracking:dyn:lin}
\begin{align}
\Delta \dot x &= A(t) \Delta x+  B(t)\Delta u + O(\|\Delta x\|^2,  \|\Delta u\|^2)\label{tracking:dyn:lin:a} \\
\Delta y &=  C(t) \Delta x + O(\|\Delta x\|^2)\label{tracking:dyn:lin:b}
\end{align}
\end{subequations}
where 
\begin{align*}
A(t) &= \frac{\partial \tilde X}{\partial x}(x_0(t),u_0(t)), \,\,B(t) = \frac{\partial \tilde X}{\partial u}(x_0(t),u_0(t)),\\
C(t) &= \frac{\partial h}{\partial x}(x_0(t)).
\end{align*}
It is understood that the big Oh $O (\cdot )$ also depends explicitly on  $t$. The following theorem is a simple application of the Lyapunov linearization method.
\begin{theorem}\label{thm:tracking:review}
If there is a time-varying gain $K(t)$ such that the LTV system
\begin{equation}\label{LTV:z}
\dot z = (A(t) - B(t)K(t))z
\end{equation}
is exponentially stable, then  the linear controller $\Delta u_0 = -K(t) \Delta x$ makes $\Delta x = 0$ an exponentially stable equilibrium point for the nonlinear closed-loop tracking error dynamics \eqref{tracking:dyn}. In consequence, the linear controller 
\begin{equation}\label{stable:LTV:controller}
u = u_0(t) - K(t) (x- x_0(t))
\end{equation}
 enables the system $\tilde \Sigma_{{\mathbb R}^n}$ in \eqref{Sigma:tilde} to  exponentially track the reference $x_0(t)$ for any initial state $x(0)$ in a neighborhood of $x_0(0)$ in ${\mathbb R}^n$.  Furthermore, the same controller  \eqref{stable:LTV:controller}, if restricted to $M$, enables the system $\Sigma_M$ to exponentially track the reference $x_0(t)$ for any initial state $x(0)$ in a neighborhood of $x_0(0)$ in $M$.
\end{theorem}

%
%
%

\subsection{Observer-Based Tracking Controllers}
We build observer-based tracking controllers for the system $\Sigma_M$ for the reference $(x_0(t),u_0(t))$ by designing observer-based tracking controllers for the ambient system $\tilde \Sigma_{{\mathbb R}^n}$ in ${\mathbb R}^n$ for the same reference trajectory and then restricting them to $M$.

\subsubsection{Linear Observer-Based Tracking Controllers}
We build a linear observer for the tracking error $\Delta x$.  Consider the following observer for the tracking error dynamics \eqref{tracking:dyn} or \eqref{tracking:dyn:lin}:
\begin{equation}\label{linear:obs:M}
\dot z_{\rm o} = A(t) z_{\rm o}+  B(t)\Delta u - L(t) (C(t)z_{\rm o}- \Delta y),
\end{equation}
where $z_{\rm o}\in {\mathbb R}^n$ is the estimate of the tracking error $\Delta x$;  $\Delta u$ and $\Delta y$ are defined in \eqref{def:Delta:u} and \eqref{def:Delta:y}; and the observer gain $L(t)$ is given by 
\begin{equation}\label{gain:L}
L(t) = P(t)C^T(t)R^{-1}(t),
\end{equation}
where $P(t)$ is the solution to 
\begin{align}
\dot P &= PA^T(t) + A(t) P  - PC^T(t)R^{-1}(t) C(t)P + Q(t), \label{Riccati:P:obs}
\end{align}
where $R(t) =R^T(t)>0$ and $Q(t) =Q^T(t)\geq 0$ are to be chosen. Let
\[
e_{\rm o} = \Delta x - z_{\rm o},
\]
which satisfies
\begin{equation}\label{eo:dynamics}
\dot e_{\rm o} = (A(t) - L(t)C(t)) e_{\rm o} + O(\|\Delta x\|^2, \|e_{\rm o}\|^2, \|\Delta u\|^2).
\end{equation}
Apply  to the system $\tilde \Sigma_{{\mathbb R}^n}$ a controller of the form 
\begin{equation}\label{obc}
u = u_0(t) - K(t) z_{\rm o}.
\end{equation}
Then the tracking error dynamics \eqref{tracking:dyn:lin} and the observation error dynamics \eqref{eo:dynamics} can be written together as
\begin{align}
\begin{bmatrix}
\Delta \dot x \\
\dot e_{\rm o}
\end{bmatrix} &= \begin{bmatrix}
A(t) - B(t)K(t) & B(t)K(t) \\
0 & A(t) - L(t)C(t)
\end{bmatrix}
\begin{bmatrix}
\Delta x \\
e_{\rm o}
\end{bmatrix} \nonumber \\ 
&\quad + O(\|\Delta x\|^2, \|e_{\rm o}\|^2).\label{total:error:dynamics}
\end{align}
\begin{theorem}\label{thm:main:obs}
Suppose that $K(t)$ is chosen such that the LTV system \eqref{LTV:z} is exponentially stable, that $B(t)K(t)$ is bounded, and that all the hypotheses  in Theorem \ref{thm:observer} in the Appendix hold true. Then, $(\Delta x, e_{\rm o}) = (0,0)$ is an exponentially stable equilibrium point for \eqref{total:error:dynamics}.  In consequence, the trajectory $x(t)$ of the system $\Sigma_M$ on $M$ exponentially tracks the reference $x_0(t)$ with the observer-based controller that consists of  \eqref{obc} and \eqref{linear:obs:M} -- \eqref{Riccati:P:obs}.

\begin{proof}
{\rm It follows from Theorem \ref{thm:observer} in the Appendix and Theorem \ref{thm:tracking:review} above.}
\end{proof}
\end{theorem}

\begin{corollary}\label{cor:main:obs}
Suppose that $K(t)$ is chosen such that the LTV system \eqref{LTV:z} is exponentially stable, that $B(t)K(t)$ is bounded, and that all the hypotheses in Corollary \ref{cor:observer} in the Appendix hold true. Then, the same conclusion as that in Theorem  \ref{thm:main:obs} holds.

\begin{proof}
{\rm It follows from Corollary \ref{cor:observer} in Appendix and Theorem \ref{thm:tracking:review} above.}
\end{proof}
 \end{corollary}

\begin{remark}{\rm 
1. Although  Kalman's theory is here used to design the observer gain $L(t)$, one can alternatively use any other method to produce $L(t)$ as far as the system $\dot e_{\rm o} = (A(t) - L(t)C(t)) e_{\rm o} $ is exponentially stable. For the rigid body system, we will build an observer of a non-Kalman type.

2. It is easy to verify that the signal defined by $
x_{\rm est}(t) = x_0(t) + z_{\rm o}(t)$
 converges to $x(t)$ as $t$ tends to infinity. In this sense, \eqref{linear:obs:M} can be regarded as a state observer for $\tilde \Sigma$.
}\end{remark}

\subsubsection{Nonlinear Observer-Based Tracking Controllers}
Consider the following state observer  for $\tilde \Sigma_{{\mathbb R}^n}$:
\begin{equation}\label{non:observer}
\dot {\hat x} = \tilde X(\hat x, u) - L(t) (h(\hat x) - y),
\end{equation}
where $\hat x$ denotes the state estimate of $x$, and   $y$ is the output of $\tilde \Sigma_{{\mathbb R}^n}$. The observer gain $L(t)$ in \eqref{non:observer} is given in \eqref{gain:L} and obtained via \eqref{Riccati:P:obs}. 
  Let
\[
e_{\rm o} = x - \hat x 
\]
denote the  error of observation of $x$. It is straightforward to show that it satisfies \eqref{eo:dynamics}.  Apply a controller of the form 
\begin{equation}\label{obc:non}
u = u_0(t) - K(t) (\hat x - x_0(t))
\end{equation}
to the system $\tilde \Sigma_{{\mathbb R}^n}$. Then the tracking error dynamics \eqref{tracking:dyn:lin} and the observation error dynamics \eqref{eo:dynamics} satisfy \eqref{total:error:dynamics}. 
\begin{theorem}\label{thm:main:obs:2}
Suppose that $K(t)$ is chosen such that the LTV system \eqref{LTV:z} is exponentially stable, that $B(t)K(t)$ is bounded, and that all the hypotheses  in Theorem \ref{thm:observer} in the Appendix hold true. Then, $(\Delta x, e_{\rm o}) = (0,0)$ is an exponentially stable equilibrium point for \eqref{total:error:dynamics}.  In consequence, the trajectory $x(t)$ of the system $ \Sigma_M$ on $M$ exponentially tracks the reference $x_0(t)$ with the observer-based controller that consists of  \eqref{obc:non}, \eqref{non:observer} and \eqref{gain:L} -- \eqref{Riccati:P:obs}.
\end{theorem}

\begin{corollary}\label{cor:main:obs:2}
Suppose that $K(t)$ is chosen such that the LTV system \eqref{LTV:z} is exponentially stable, that $B(t)K(t)$ is bounded, and that all the hypotheses in Corollary \ref{cor:observer} in the Appendix hold true, Then, the same conclusion as that in Theorem  \ref{thm:main:obs:2} holds. 
 \end{corollary}

\begin{remark}{\rm
1.  The difference between \eqref{linear:obs:M} and \eqref{non:observer} is that \eqref{linear:obs:M} estimates the tracking error $\Delta x = x - x_0(t)$ whereas  \eqref{non:observer} estimates the state $x$.   However, the underlying error dynamics share the same first-order approximation as shown in  \eqref{total:error:dynamics}.

2.  It is noteworthy that the {\it high-gain} observer  \eqref{non:observer} is here linearized along the reference trajectory $x_0(t)$.  Since the reference $x_0(t)$ is chosen first, it is possible to verify the conditions in Theorems \ref{thm:main:obs} and \ref{thm:main:obs:2} and Corollaries \ref{cor:main:obs} and \ref{cor:main:obs:2} in advance, which also allows to design the time-varying observer gain $L(t)$ in advance by integrating \eqref{Riccati:P:obs}. 
}\end{remark}
\section{Application: The Fully Actuated Rigid Body System}

\subsection{Review of  the Embedding, Transversal Stabilization and Linearization of the Rigid Body System}
We  review the embedding of the rigid body system from $\SO \times {\mathbb R}^3$ into ${\mathbb R}^{3\times 3} \times {\mathbb R}^3$; refer to \cite{Ch17} for more detail.   Let $\SO = \{ R\in {\mathbb R}^{3\times 3} \mid R^TR = I, \det R >0 \}$ be the set of all $3 \times 3$ rotation matrices and ${\mathfrak so}(3) = \{A \in {\mathbb R}^{3\times 3} \mid A = -A^T\} $ the set of all $3\times 3$ skew symmetric matrices.  The hat map $\hat{} : {\mathbb R}^3 \rightarrow \mathfrak{so}(3)$ is defined as follows: 
\[
\hat \Omega = \begin{bmatrix}
0 & -\Omega_3 &\Omega_2 \\
\Omega_3 & 0 & -\Omega_1\\
-\Omega_2 & \Omega_1 & 0
\end{bmatrix}
\]
for  all $\Omega = (\Omega_1, \Omega_2, \Omega_3) \in {\mathbb R}^3$.
The hat map satisfies the identity, $\hat x y = x \times y$ for all $x,y \in {\mathbb R}^3$.  The inverse map of the hat map is called the vee map and denoted by $\vee$ so that $(\hat{\Omega})^\vee = \Omega$ for all $\Omega \in {\mathbb R}^3$.  We use the inner product $\langle A, B \rangle = \operatorname{tr}(A^TB)$ for $A,B\in {\mathbb R}^{n\times n}$, and $\| \cdot \|$ denotes the norm induced from this inner product on ${\mathbb R}^{n\times n}$. The symbol  $[\, , ]$ denotes the usual matrix commutator: $[A,B] = AB - BA$  for all $A,B\in {\mathbb R}^{n\times n}$. The operators $\Sym$ and $\Skew$ denote the symmetrization operator and the skew-symmetrization operator, respectively, on square matrices, i.e. 
\[
\Sym (A) = \frac{1}{2}(A+A^T), \quad \Skew (A) = \frac{1}{2}(A-A^T)
\]
for any square matrix $A$.

The equations of motion of the fully actuated rigid body system are given by
\begin{subequations}\label{rigid:original}
\begin{align}
\dot R &= R\hat \Omega, \\
\dot \Omega &= \MI^{-1} ( \MI \Omega \times \Omega) + \MI^{-1} u,\\
y &= R,
\end{align}
\end{subequations}
where $(R,\Omega) \in \SO \times  {\mathbb R}^3 \subset {\mathbb R}^{3\times 3} \times {\mathbb R}^3$ is the state vector consisting of a rotation matrix $R$ and a body angular velocity $\Omega$; $u \in {\mathbb R}^3$ is the control torque;  $\mathbb I$ is the moment of inertial matrix of the rigid body; and $y$ is the output of the system.  The above dynamics naturally extend to ${\mathbb R}^{3\times 3} \times {\mathbb R}^3$ as they are, by considering $R$ as a $3\times 3$ matrix. 
Consider a   function $\tilde V$ on $\operatorname{GL}^+(3) \times {\mathbb R}^3\subset { \mathbb R}^{3\times 3} \times {\mathbb R}^3$ that is defined by 
\[
\tilde V(R,\Omega) = \frac{\coV }{4}\|R^TR - I\|^2,
\]
where $\coV >0$ and $\operatorname{GL}^+(3) = \{ R\in {\mathbb R}^{3\times 3} \mid \det R >0\}$. It satisfies  $\tilde V^{-1}(0) = \SO\times {\mathbb R}^3$  and \eqref{tilde:V:X} for all $(R,\Omega,u) \in \operatorname{GL}^+(3)  \times {\mathbb R}^3 \times {\mathbb R}^3$.  The gradient $\nabla \tilde V = (\nabla_R \tilde V, \nabla_\Omega \tilde V)$ is computed as
\[
\nabla_R \tilde V = -\coV R(R^TR - I), \quad \nabla_\Omega \tilde V = 0.
\]
With this function $\tilde V$,  the system corresponding to \eqref{Sigma:tilde} is 
\begin{subequations}\label{rigid:tilde:eq}
\begin{align}
\dot R &= R\hat \Omega - \coV R(R^TR - I), \label{R:s:eq}\\
\dot \Omega &= \MI^{-1} ( \MI \Omega \times \Omega) + \MI^{-1} u,\\
y&=R,
\label{Omega:s:eq}
\end{align}
\end{subequations}
where $(R,\Omega) \in {\mathbb R}^{3\times 3} \times { \mathbb R}^3$. It is straightforward to show that Theorem \ref{theorem:attractive:M} holds for \eqref{rigid:tilde:eq}, so   $\SO \times {\mathbb R}^3$ is an exponentially stable invariant manifold of \eqref{rigid:tilde:eq}. It is trivial to see that the system \eqref{rigid:tilde:eq} reduces to \eqref{rigid:original} on $\SO \times {\mathbb R}^3$.

Take a reference trajectory $(R_0(t), \Omega_0 (t)) \in \SO \times { \mathbb R}^3$ and the corresponding control signal $u_0(t)$ such that 
\begin{subequations}
\begin{align}
\dot R_0(t) &= R_0(t)\hat \Omega_0 (t), \label{ref:traj:rigid:a} \\
\dot \Omega_0 (t) &= \MI^{-1} ( \MI \Omega_0(t) \times \Omega_0(t)) + \MI^{-1} u_0 (t) \label{ref:traj:rigid:b}\\
y_0(t) &:= R_0(t)
\end{align}
\end{subequations}
for all $t\geq 0$. Assume that  $(R_0(t), \Omega_0(t))$ and $u_0(t)$ are bounded  over the time interval $[0,\infty )$. The paper \cite{Ch17} provides several tracking controllers for this type of reference trajectories.
Let
\begin{align*}
&\Delta R(t)= R(t)-R_0(t) \in {\mathbb R}^{3\times 3}, \\
&\Delta \Omega(t)= \Omega(t)-\Omega_0(t) \in { \mathbb R}^3,\\
&\Delta u(t)=u(t)-u_0(t) \in{\mathbb R}^3,\\
&\Delta y = y - y_0(t) = R - R_0(t)  \in {\mathbb R}^{3\times 3}
\end{align*}
denote tracking errors. Then, the tracking error dynamics can be written as
\begin{subequations}\label{rigid:tracking:error:dyn}
{\small
\begin{align}
\Delta \dot R  &= \Delta R \hat \Omega_0  + R_0  \Delta \hat \Omega -2 \coV R_0  \text{Sym} (R_0^T \Delta R) + O(2), \\
\Delta \dot \Omega  &= \MI^{-1} ( \MI \Delta \Omega \times \Omega_0+\MI  \Omega_0 \times  \Delta \Omega) + \MI^{-1} \Delta u + O(2),\\
\Delta y &= \Delta R,
\end{align}}
\end{subequations}
where  $O(2)=O(\|\Delta R\|^2, \|\Delta \Omega\|^2, \|\Delta u\|^2)$. 
Introduce a new matrix variable $Z$ to replace $\Delta R$ as follows:
\[
Z = R_0(t)^T \Delta R.
\]
 Let
\[
Z_s = \Sym (Z), \quad Z_k=\Skew (Z)
\]
such that $Z = Z_s + Z_k$.
Then,  the tracking error dynamics are transformed to
\begin{subequations}\label{Zs:Zk:Om:tv}
{\small
\begin{align}
\dot Z_s &= [Z_s, \hat \Omega_0] - 2\coV Z_s + O(2),\label{Zs:Zk:Om:tv:a}\\
\dot Z_k^\vee &= Z_k^\vee \times \Omega_0 + \Delta \Omega +O(2), \label{Zs:Zk:Om:tv:b}\\
\Delta \dot \Omega &=  \MI^{-1} ( \MI \Delta \Omega \times \Omega_0+\MI  \Omega_0 \times  \Delta \Omega) + \MI^{-1} \Delta u   +O(2),
 \label{Zs:Zk:Om:tv:c}\\
 \Delta y &= R_0(t)(Z_s + Z_k),
\end{align}
}
\end{subequations}
where $
O(2)=O(\|Z\|^2, \|\Delta \Omega\|^2, \|\Delta u\|^2). $
Since the reference $R_0(t)$ is known, the output $\Delta y$ can be replaced with $(\Delta y_s, \Delta y_k)$ that are defined by
\begin{subequations}\label{output:sk}
\begin{align}
\Delta y_s  &= \Sym(R_0(t)^T\Delta y) = Z_s \in \Sym ({\mathbb R}^{3\times 3}),\label{output:sk:s}\\
\Delta y_k &= \Skew(R_0(t)^T\Delta y)^\vee = Z_k^\vee \in {\mathbb R}^3. \label{output:sk:k}
\end{align}
\end{subequations}

\subsection{Observer-based Tracking Controller Design}

We choose to use the linear tracking error observer \eqref{linear:obs:M} to build an observer-based tracking controller for the rigid body system with the measurement of $R$.  The linear part of \eqref{Zs:Zk:Om:tv:a} is already exponentially stable and decoupled from the rest of the dynamics, so there is no need to stabilize it. So, we have only to focus on exponentially stabilizing the linear part of \eqref{Zs:Zk:Om:tv:b} and \eqref{Zs:Zk:Om:tv:c}. Hence, it suffices to build an observer for $(  Z_k^\vee, \Delta \Omega)$ with the output $\Delta y_k$.  In view of the linear part of \eqref{Zs:Zk:Om:tv:b} and \eqref{Zs:Zk:Om:tv:c}, the linear observer corresponding to  \eqref{linear:obs:M} is written as
\begin{equation}\label{full:linearized:rigid}
\dot z_{\rm o} = A(t) z_{\rm o} + B \Delta u + L(t) (C z_{\rm o} - Z_k^\vee), 
\end{equation}
where 
\[
z_{\rm o} = (Z^\vee_{k,\rm est}, \Delta \Omega_{\rm est}) \in {\mathbb R}^3 \times { \mathbb R}^3
\]
 is the estimate of  $(Z_k^\vee, \Delta \Omega) \in { \mathbb R}^3 \times {\mathbb R}^3$, 
and
\begin{align*}
A(t) &= \begin{bmatrix}
-\hat \Omega_0(t) &  I \\
0 & {\mathbb I}^{-1}( \widehat{{\mathbb I} \Omega_0(t)} -\hat\Omega_0(t) { \mathbb I })
\end{bmatrix}\!\!, B = \begin{bmatrix} 
0 \\ {\mathbb I}^{-1}
\end{bmatrix},\\
C &= \begin{bmatrix} I & 0 \end{bmatrix}.
\end{align*}

\begin{lemma}\label{lem:ucb}
If $\Omega_0(t)$ is periodic, then $(A(t), C)$ is uniformly completely observable.  
%
\end{lemma}

\begin{lemma}\label{lem:ucb:2}
If $\Omega_0(t)$ is periodic, then the pair $(A(t),I_{6\times 6})$ is uniformly completely controllable.
\end{lemma}

By Theorem II.5 in \cite{Ch17}, any  controller of the form
\begin{align*}
u &= u_0- ({\mathbb I} \Delta \Omega)\times \Omega_0 - ({\mathbb I }\Omega_0) \times \Delta \Omega \nonumber \\&\quad - {\mathbb I }(k_P Z_k^\vee +K_D\Delta \Omega)
\end{align*}
with $k_P>0$ and $K_D = K_D^T \in {\mathbb R}^{3\times 3}$ positive definite,
exponentially stabilizes the tracking error dynamics \eqref{Zs:Zk:Om:tv}. This form of controller leads to the following observer-based tracking controller:
\begin{align}\label{ob:tracking:rigid}
u &= u_0- ({\mathbb I} \,\Delta \Omega_{\rm est})\times \Omega_0 - ({\mathbb I } \Omega_0) \times\Delta \Omega_{\rm est} \nonumber \\
&\qquad - { \mathbb I } (k_P Z^\vee_{k,\rm est} +K_D\Delta \Omega_{\rm est})
\end{align}
where $(Z^\vee_{k,\rm est}, \Delta \Omega_{\rm est}) = z_{\rm o} $ is obtained from \eqref{full:linearized:rigid} and 
\[
L(t)  = P(t)C^T {\tilde R}^{-1}(t)
\]
where $P(t)$ is the solution to
\begin{align}\label{rigid:observer:Riccati}
\dot P &= PA^T(t) + A(t) P  - PC^T\tilde R^{-1}(t) CP(t) + Q(t),
\end{align}
where $\tilde R(t) =\tilde R^T(t)>0$ and $Q(t) = Q^T(t)\geq 0$ are chosen such that there are positive numbers $\gamma_i$, $i =1,\ldots, 4$ such that 
\[
\gamma_1 I \leq Q(t) \leq \gamma_2 I, \quad \gamma_3 I \leq \tilde R(t) \leq \gamma_4 I 
\]
for all $t\geq 0$. Here, we intentionally put a tilde over $R$ in the above three equations since $R$ is reserved for rotation matrix in this section. From  Theorem \ref{thm:main:obs} and Lemmas \ref{lem:ucb} and \ref{lem:ucb:2}, we obtain the following theorem:
\begin{theorem}
The  controller \eqref{ob:tracking:rigid} exponentially stabilizes the tracking error dynamics \eqref{Zs:Zk:Om:tv} or \eqref{rigid:tracking:error:dyn} if $\Omega_0(t)$ is periodic.
\end{theorem}
\begin{remark}
Notice that we only need
\[
\Delta y_k = \Skew (R_0(t)^T\Delta R)^\vee =  \Skew (R_0(t)^T R)^\vee 
\]
for the observer designed above, instead of   full information on the rotation matrix $R$. 
\end{remark}

We now build a non-Kalman type observer gain $L(t)$   which does not require the periodicity of $\Omega_0(t)$.  
\begin{lemma}\label{lemma:nonKalman:obs}
Let 
\[
L(t) = \begin{bmatrix} L_1(t) \\ L_2(t) \end{bmatrix}, 
\]
where
\begin{align*}
L_1(t) &= -\hat\Omega_0(t) + {\mathbb I}^{-1}( \widehat{{\mathbb I} \Omega_0(t)} -\hat\Omega_0(t) { \mathbb I }) + M_1,\\
L_2(t) &= -({\hat{\dot \Omega}}_0(t) + \dot L_1(t)) + (\hat \Omega_0(t) + L_1(t))^2 \\
&\quad - M_1 (\hat \Omega_0(t) + L_1(t)) + M_2
\end{align*}
with   any constant $3\times 3$ positive definite symmetric matrices $M_1$ and $M_2$.  Then, the observation error dynamics $
\dot e_{\rm o} = (A(t) - L(t)C) e_{\rm o}$
for the observer \eqref{full:linearized:rigid} is exponentially stable.  Here, $e_{\rm o} = (Z_k^\vee, \Delta \Omega) - (Z^\vee_{k,\rm est}, \Delta \Omega_{\rm est})$.
\end{lemma}
\begin{theorem}
The  controller \eqref{ob:tracking:rigid} with the observer gain provided in Lemma \ref{lemma:nonKalman:obs} exponentially stabilizes the tracking error dynamics \eqref{Zs:Zk:Om:tv} or \eqref{rigid:tracking:error:dyn}.
\end{theorem}
\subsection{Simulation}
\begin{figure}[tb]
\vspace{-3mm}
\begin{center}
\includegraphics[scale = 0.23]{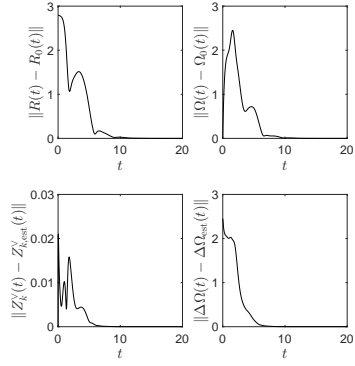}
\end{center}
\vspace{-3mm}
\caption{\label{figure.tracking:rigid} The tracking and observation errors.  }
\end{figure}

The moment of inertia matrix of the system is given by $
{\mathbb I } = \operatorname{diag}[3,2,1]$.
The parameter $\coV$ in \eqref{R:s:eq} is set to 1, and 
the control parameters in \eqref{ob:tracking:rigid} are chosen as $k_P=4$ and  $ K_D = 4I$.
The reference trajectory $(R_0(t), \Omega_0(t)) \in \SO \times {\mathbb R}^3$ with the reference control signal $u_0(t) \in { \mathbb R}^3$ are 
\begin{align*}
R_0(t)  &\\
&\!\!\! \!\!\!\!\!\! \!\!\!\!{\footnotesize = \begin{bmatrix}
\cos^2 t & - \sin t & \cos t\sin t \\
\sin^2 t + \cos^2 t\sin t & \cos^2 t & \cos t \sin^2 t - \cos t \sin t \\
\cos t \sin^2t - \cos t \sin t & \cos t \sin t & \cos^2t + \sin^3t 
\end{bmatrix}},\\
\Omega_0 (t)&= [
\cos^2 t - \sin t,  1-\sin t,  (1+  \sin t) \cos t)]^T,\\
 u_0(t) &= {\mathbb I} \dot \Omega_0 - {\mathbb I} \Omega_0 \times \Omega_0=
{\small  \begin{bmatrix}
-(3+6\sin t + \cos^2t )\cos t \\ -2(2+\sin t)\cos t \sin^2t \\ -(2\sin t  - \cos^2 t)\sin t
\end{bmatrix}. }
\end{align*}
which satisfy \eqref{ref:traj:rigid:a} and  \eqref{ref:traj:rigid:b}. The initial condition is given by $
R(0) = \exp (0.9\pi \hat e_2)$, and  $\Omega (0) =  (1,1,1)$,
where $R(0)$ is a rotation around $e_2 = (0,1,0)$ through $0.9\pi$ radians. We then have the initial attitude tracking error $\|R(0) - R_0(0)\| = 2.7936$ which is  fairly close to $2\sqrt 2 = 2.8284$, the magnitude of  maximum possible tracking error. The initial state for  the tracking error observer is set to $
Z_{k, \rm est}^\vee (0) =(0,0,0)$ and $ \Delta \Omega_{\rm est} (0) = (1,2,1)$. 
We choose the following values of observer parameters in \eqref{rigid:observer:Riccati}: $
 Q (t)\equiv 100 I$, $\tilde R (t) \equiv 0.01I$,  $P(0) = 100I$.
The simulation results are plotted in Figure \ref{figure.tracking:rigid}, where it can be seen that both the tracking error and the observation error converge to zero as time tends to infinity.

\appendix
We review the theory of linear time varying systems, following \cite{Ch84,Ch99,Ka60,KaBu61} for easy reference.
Consider a linear time varying (LTV) system:
\begin{align*}
\dot x = A(t)x + B(t)u,\qquad y = C(t)x,
\end{align*}
where $x\in { \mathbb R}^n$, $u\in {\mathbb R}^\dimu$, $y \in {\mathbb R}^p$; and $A(t)$, $B(t)$ and $C(t)$ are continuously time-varying matrices of appropriate dimensions.  Let $\Phi (t,\tau) \in { \mathbb R}^{n\times n}$ denote the transition matrix of the system with $u = 0$, i.e. the solution to the initial value problem
\[
\frac{\partial \Phi}{\partial t} (t,\tau)  = A(t)\Phi(t,\tau), \quad \Phi(\tau,\tau) = I.
\]

\begin{definition}\label{def:ucc}
The LTV system, or the pair $(A(t), B(t))$, is  uniformly completely controllable if there is a positive number $\sigma$ and positive numbers $\alpha_i$, $i = 1, \ldots, 4$ that depend on $\sigma$ such that 
\begin{align}
&\alpha_1 I \leq W(t,t+\sigma) \leq \alpha_2 I, \label{ucc1}\\
& \alpha_3  I \leq \Phi(t+\sigma, t) W(t,t+\sigma)\Phi^T(t+\sigma, t)  \leq \alpha_4I  \label{ucc2}
\end{align}
for all $t$,
where  the matrix $W(\cdot, \cdot)$ is defined as 
\[
W(t,\bar t) = \int_t^{\bar t} \Phi (t,\tau) B(\tau) B^T(\tau) \Phi^T(t,\tau) d\tau.
\]
for $t,\bar t \in \mathbb R$.
\end{definition}

\begin{definition}\label{def:uco}
The LTV system or the pair of matrices $(A(t), C(t))$, is  uniformly completely observable if there is a positive number $\sigma$ and positive numbers $\alpha_i$, $i = 1, \ldots, 4$ that depend on $\sigma$ such that 
\begin{align}
&\alpha_1  I \leq V(t,t+\sigma) \leq \alpha_2 I, \label{uco1}\\
& \alpha_3  I \leq \Phi^T(t,t+\sigma) V(t,t+\sigma)\Phi(t,t+\sigma) \leq \alpha_4I \label{uco2}
\end{align}
for all $t$,
where  the matrix $V(\cdot, \cdot)$ is defined as 
\begin{equation}\label{def:V}
V(t,\bar t) = \int_t^{\bar t} \Phi^T (\tau, t) C^T(\tau) C(\tau) \Phi(\tau, t) d\tau.
\end{equation}
for $t,\bar t \in \mathbb R$.
\end{definition}

\begin{lemma}\label{lem:ucc:uco}
Suppose that $A(t)$ is bounded. Then the following hold:

1. If one of  \eqref{ucc1} and \eqref{ucc2}   in Definition \ref{def:ucc} holds, then the other hold too. In other words, if there exist $\alpha_1$ and $\alpha_2$ such that \eqref{ucc1} holds, then there exist $\alpha_3$ and $\alpha_4$ such that \eqref{ucc2} holds; and vice versa.

2. If one of \eqref{uco1} and \eqref{uco2}   in Definition \ref{def:uco} holds, then the other holds too.

\begin{proof}
See \cite{Ka60}.
\end{proof}
\end{lemma}

Consider the following observer system:
\begin{align*}
\dot {\hat x} &= A(t) \hat x + B(t) u - L(t)(C(t) \hat x - y),\\
\dot P &= PA^T(t) + A(t) P - PC^T(t)R^{-1}(t) C(t)P(t) + Q(t),
\end{align*}
where $R(t) =R^T(t)>0$, $Q(t) =Q^T(t)\geq 0$, and 
\[
L(t) = P(t)C^T(t)R^{-1}(t).
\]
Then, the observation error $e_{\rm o} = x(t)  - \hat x(t)$ satisfies
\begin{equation}\label{eo:dynamics:app}
\dot e_{\rm o} = (A(t) - L(t)C(t))e_{\rm o}. 
\end{equation}

\begin{theorem}\label{thm:observer}
Suppose that the pair $(A(t), C(t))$ is uniformly completely observable; that for some matrices $D(t)$ and $\tilde Q(t) =\tilde Q^T(t)$ such that $
Q(t) = D(t)\tilde Q(t)D^T(t)$,
 the pair $(A(t), D(t))$ is uniformly completely controllable and there are positive numbers $\gamma_1$ and $\gamma_2$ such that $ \gamma_1 I \leq \tilde Q(t) \leq \gamma_2 I$ for all $t$; and that there are positive numbers $\gamma_3$ and $\gamma_4$ such that  and $\gamma_3 I \leq  R(t) \leq \gamma_4 I$ for all $t$. Then, the observation error dynamics \eqref{eo:dynamics:app} is exponentially stable. 
\end{theorem}

The hypotheses of uniformly complete observability and controllability in the above theorem can be relaxed if boundedness of $A(t)$ is assumed. The following corollary follows from the above theorem and Lemma \ref{lem:ucc:uco}.
\begin{corollary}\label{cor:observer}
Suppose that $A(t)$ is bounded; that the pair $(A(t), C(t))$ satisfies \eqref{uco1} in Definition \ref{def:uco}; that for some matrices $D(t)$ and $\tilde Q(t) =\tilde Q^T(t)$  such that  $
Q(t) = D(t)\tilde Q(t)D^T(t)$,
 there are positive numbers $\gamma_1$ and $\gamma_2$ such that $ \gamma_1 I \leq \tilde Q(t) \leq \gamma_2 I$ for all $t$ and the pair $(A(t), D(t))$ satisfies \eqref{ucc1}, where $D(t)$ is used in place of $B(t)$ in the computation of $W(\cdot, \cdot)$; and that there are positive numbers $\gamma_3$ and $\gamma_4$ such that  and $\gamma_3 I \leq  R(t) \leq \gamma_4 I$ for all $t$. Then, the observation error dynamics \eqref{eo:dynamics:app} is exponentially stable. 
\end{corollary}

Suppose that $u = -K(t) x $ exponentially stabilizes the LTV system. If we use the state estimate $\hat x(t)$ from the observer and apply $u = - K(t) \hat x$ to the LTV system instead, then the dynamics of resulting system, which is comprised of the LTV system and the observer, is written as
\begin{equation}\label{composite:system}
\begin{bmatrix}
\dot x \\ \dot e_{\rm o}
\end{bmatrix} = \begin{bmatrix}
A(t) - B(t)K(t) & B(t)K(t) \\
0 &  A(t) - L(t)C(t)
\end{bmatrix}
\begin{bmatrix}
 x \\ e_{\rm o}
\end{bmatrix}
\end{equation}
where $e_{\rm o}= x(t)  - \hat x(t)$ is the observation error. The following theorem is elementary.
\begin{theorem}
Suppose that the system 
\[
\dot x = (A(t) - B(t)K(t))x
\]
 is exponentially stable; that the hypotheses in Theorem \ref{thm:observer} or Corollary \ref{cor:observer} hold; and that $B(t)K(t)$ is bounded. Then, the composite system \eqref{composite:system} is exponentially stable.
\end{theorem}

\subsection*{Acknowledgement}
This research has been in part supported by KAIST  under grant G04170001 and  by the ICT R\&D program of MSIP/IITP [2016-0-00563, Research on Adaptive Machine Learning Technology Development for Intelligent Autonomous Digital Companion].

%
%


\begin{thebibliography}{99}

\bibitem{Ch11}
D.E. Chang, ``A simple proof of the Pontryagin maximum principle on manifolds,'' {\it Automatica,}  47 (3), 630 -- 633, 2011.

\bibitem{Ch17} D.E. Chang, ``Controller design for systems on manifolds in {E}uclidean space,'' in {\it Proc. IEEE Conference on Decision and Control}, Melbourne, Australia, 2017. {\it arXiv preprint arXiv:1710.02780}.

\bibitem{ChJiPe16}
D.E. Chang, F. Jim\'{e}nez and M. Perlmutter, ``Feedback integrators,'' {\it J. Nonlinear Science,}  26(6), 1693 -- 1721, 2016.


\bibitem{Ch84} C.-T. Chen, {\it Linear System Theory and Design}, 2nd Ed., Oxford University Press,  New York, 1984.

\bibitem{Ch99} C.-T. Chen, {\it Linear System Theory and Design}, 3rd Ed., Oxford University Press,  New York, 1999.


\bibitem{Ka60} R.E. Kalman, ``Contribution to the theory of optimal control,'' {\it Bol. Soc. Mat. Mex.}, Vol. 5, 102 -- 119, 1960. 

\bibitem{KaBu61} R.E. Kalman and R.S. Bucy, ``New results in linear filtering and prediction theory,'' {\it Trans. ASME J. Basic Engineering},  83 (1), 95 -- 108, 1961.



\end{thebibliography}
\end{document}